\documentclass[a4paper,11pt]{article}
\usepackage{amssymb,amsmath}
\usepackage{amsthm,enumerate}
\usepackage{amsfonts}
\usepackage[latin1]{inputenc}
\usepackage{epsfig}
\usepackage[dvips]{color}

\newtheorem{theorem}{Theorem}[section]
\newtheorem{lemma}[theorem]{Lemma}

\newtheorem{definition}[theorem]{Definition}

\newtheorem{remark}[theorem]{Remark}

\newtheorem*{thma}{Theorem A}

\numberwithin{equation}{section}

\begin{document}

\newcommand{\diam}{\operatorname{diam}}
\newcommand{\Vol}{\operatorname{Vol}}
\newcommand{\dist}{\operatorname{dist}}
\newcommand{\length}{\operatorname{length}}
\newcommand{\inter}{\operatorname{int}}
\newcommand{\B}{\mathbb{B}}
\newcommand{\loc}{\operatorname{loc}}
\newcommand{\capacity}{\operatorname{cap}}

\long\def\symbolfootnote[#1]#2{\begingroup%
 \def\thefootnote{\fnsymbol{footnote}}\footnote[#1]{#2}\endgroup}

\title{Quasiconformal mappings and singularity of boundary distortion}
\author{Tomi Nieminen \and Ignacio Uriarte-Tuero}

\date{}
\maketitle

\begin{abstract}\symbolfootnote[0]{\emph{Mathematics Subjects Classification} (2000). Primary 30C65.}
We extend a well-known theorem by Jones and Makarov \cite{JM} on the singularity of boundary distortion of planar conformal mappings. We use a different technique to recover the previous result and, moreover, generalize the result for quasiconformal mappings of the unit ball $\B^n\subset \mathbb{R}^n$, $n\ge 2$. We also establish an estimate on the Hausdorff (gauge) dimension of the boundary of the image domain outside an exceptional set of given size on the sphere $\partial \B^n$. Furthermore, we show that this estimate is essentially sharp.
\end{abstract}

\section{Introduction}

Let $f:\mathbb{D}\to \Omega$ be a conformal mapping of the unit disk $\mathbb{D}$ onto a domain $\Omega\subset \mathbb{C}$. Recall that, by a classical theorem of Beurling \cite[p. 215]{Po}, the boundary function of $f$ is defined in terms of angular limits everywhere on $\partial \mathbb{D}$ except for a set of zero logarithmic capacity. Some time ago, Jones and Makarov \cite{JM} established the following remarkable result considering the singularity of boundary distortion of $f$. They write $f^*\Lambda_{\varphi}\perp m_2$, if $f$ maps the whole unit circle except a set of zero $\Lambda_{\varphi}$-measure onto a set of zero area. Here $\Lambda_{\varphi}$ denotes the Hausdorff measure on $\partial \mathbb{D}$ associated to a weight function $\varphi$, see below for the definition of this measure.
\begin{thma}[\cite{JM}]\label{jonestheorem}
Let $\varphi$ be a weight function satisfying $\varphi(2r)\le C\varphi(r)$, $r>0$. Then the relation $$f^*\Lambda_{\varphi}\perp m_2$$
holds for every univalent function $f$ if and only if
\begin{equation}\label{jonesdivergence}
\int_0\Big|\frac{\log \varphi(t)}{\log t}\Big|^2 \frac{dt}{t} = \infty.
\end{equation}
\end{thma}

In this note we extend the above result for quasiconformal mappings of the unit ball $\B^n$, $n\ge 2$, of Euclidean space. By definition, a homeomorphism $f:\B^n\to \Omega\subset \mathbb{R}^n$ is $K$-quasiconformal if $f\in W_{\loc}^{1,n}(\B^n;\Omega)$ and the inequality
\begin{equation}\label{distortion}
|Df(x)|^n\le KJ_f(x)
\end{equation}
holds for almost every $x\in \B^n$. Here $|Df(x)|$ stands for the operator norm of the differential matrix of $f$ at the point $x$, while $J_f(x)$ denotes the determinant of $Df(x)$. Recall that, by the analog of Beurling's theorem, the boundary mapping of $f$ is defined in terms of radial limits everywhere on $\partial \B^n$ except for a set of zero conformal ($n$)-capacity, see e.g. \cite[Theorem 4.4]{BKR}. In this setting we establish the following theorem.
\begin{theorem}\label{maincorollary}
Let $\varphi(t)$ be a weight function satisfying the technical conditions \eqref{prop1}, \eqref{prop2} and \eqref{prop3} below and denote $u=\varphi^{-1}$. Then the relation
\begin{equation}\label{claim1}
f^*\Lambda_{\varphi}\perp m_n
\end{equation}
holds for every quasiconformal mapping $f:\B^n\rightarrow \Omega\subset \mathbb{R}^n$ if and only if
\begin{equation}\label{divergence}
\int_0\Big(\frac{u(t)}{u'(t)}\Big)^{n-1}\frac{dt}{t^n}=\infty.
\end{equation}
\end{theorem}
Note that, in the case $n=2$, the condition \eqref{divergence} is equivalent to the condition \eqref{jonesdivergence}, see \cite[Remark 5.3]{N}. Thus, in the planar case, we recover the result of Jones and Makarov. Our assumptions on the weight function $\varphi$ are described more precisely in the following. We will assume that for all sufficiently small $t>0$ the function $\varphi(t)$ is an increasing and differentiable function, which satisfies $\varphi(0)=0$, $\varphi(2t)\le \beta \varphi(t)$, and
\begin{equation}\label{prop1}
\frac{\varphi'(t)t\log t}{\varphi(t)\log\varphi(t)}\ \mbox{is non-increasing or non-decreasing},
\end{equation}
and
\begin{equation}\label{prop2}
\frac{u(t)}{tu'(t)}\ \mbox{is non-decreasing},
\end{equation}
and
\begin{equation}\label{prop3}
\log\frac{1}{u(t^2)} \le \beta\log\frac{1}{u(t)}
\end{equation}
for $u=\varphi^{-1}$ with some constant $\beta>1$. Note that these technical assumptions are harmless in the sense that they are satisfied in all interesting situations, see e.g. Remark \ref{example} below. It is the condition \eqref{divergence} that is interesting in Theorem \ref{maincorollary}. Let us point out that some conditions on the regularity of $\varphi$ are assumed also in the proof of Theorem A, see \cite[p. 447-448]{JM}.

We also extend the above result by replacing the $n$-dimensional Lebesgue measure $m_n$ in \eqref{claim1} with a measure $\Lambda_{\psi}$, where the gauge function $\psi$ depends on $\varphi$. Recall that the \emph{generalized Hausdorff measure} $\Lambda_{\varphi}$ (or simply $\varphi$\emph{-measure}) is defined by
$$\Lambda_{\varphi}(E)=\lim_{r\to 0}\Big(\inf\Big\{\sum \varphi(\diam B_i): E\subset \bigcup B_i,\ \diam(B_i)\le r \Big\}\Big),$$
where the dimension gauge function $\varphi$ is required to be continuous and increasing with $\varphi(0)=0$. In particular, if $\varphi(t)=t^{\alpha}$ with some $\alpha>0$, then $\Lambda_{\varphi}$ is the usual \emph{$\alpha$-dimensional Hausdorff measure} denoted also by $H^{\alpha}$. See \cite{F} or \cite{R} for more information on the generalized Hausdorff measure.

Our main result is the following theorem. We write $f^*\Lambda_{\varphi}\perp \Lambda_{\psi}$, if $f$ maps the whole unit sphere except a set of zero $\varphi$-measure to a set of zero $\psi$-measure.
\begin{theorem}\label{main}
Let $\varphi$ be a weight function such that $u=\varphi^{-1}$ satisfies the condition \eqref{divergence} in addition to the technical properties \eqref{prop1}, \eqref{prop2} and \eqref{prop3}. Then the relation
\begin{equation}\label{claim2}
f^*\Lambda_{\varphi}\perp \Lambda_{\psi}
\end{equation}
holds for every K-quasiconformal mapping $f:\B^n\to \Omega\subset \mathbb{R}^n$ if there are positive constants $r_0$ and $C_1$ so that
\begin{equation}\label{gauge}
\psi(r)\le C_1r^n\exp\Big(C_2\int_{[r,
r_0]}\Big(\frac{u(t)}{u'(t)}\Big)^{n-1}\frac{dt}{t^n}\Big)
\end{equation}
for all $r<r_0$. Here $C_2=C_2(n,K,\beta)>0$.
\end{theorem}

We will also show in Section \ref{sharpnesssection} that Theorem \ref{main} is sharp in the following sense. Suppose that $\varphi$ satisfies the conditions of Theorem \ref{main}. Then there is an open set $\Omega$, a constant $\tilde{C}>0$, and a quasiconformal mapping $f:\B^n\to \Omega$ so that, for any set $E\subset \partial \B^n$ with $\Lambda_{\varphi}(E)=0$, we have $\Lambda_{\psi}(f(\partial \B^n\setminus E))>0$ with a dimension gauge $\psi$ satisfying
$$\psi(r)\le r^n\exp\Big(\tilde{C}\int_{[r,r_0]}\Big(\frac{u(t)}{u'(t)}\Big)^{n-1}\frac{dt}{t^n}\Big)$$
for all sufficiently small $r>0$, provided that $\Lambda_{\psi}$ is absolutely continuous with respect to $H^1$, which is the interesting case.

Let us close this section with a concrete example of our results. This remark also demonstrates which dimension gauge functions are critical for the condition \eqref{claim1} to hold.
\begin{remark}\label{example}
Let $s\ge 1$ and let
$$\varphi(t)=\exp\Big(-\Big(c\log\frac{1}{t}\Big)^{1/s}\Big).$$
Let $f$ be a $K$-quasiconformal mapping of the unit ball $\B^n$, $n\ge 2$. Then the condition \eqref{claim1} holds if and only if $s\le n/(n-1)$. In the case $s=n/(n-1)$ the condition \eqref{claim2} holds with a dimension gauge function
\begin{equation}\label{concretepsi}
\psi(r)=r^n\Big(\log\frac{1}{r}\Big)^C,
\end{equation}
where $C>0$ depends only on $n,K$ and $c$. Note that the latter statement is stronger than \eqref{claim1}.
\end{remark}

\section{Proof of the main result}

In this section we prove Theorem \ref{main} and the ``if''-part of Theorem \ref{maincorollary} (observe that the latter immediately follows from Theorem \ref{main}). Most of the machinery needed for the proof has already been developed in \cite{N}, including the following geometric concept which we use as a tool.
\begin{definition}
Let $E\subset \mathbb{R}^n$ be a compact set. Let $\alpha:(0,1)\to (0,1)$ be a continuous function such that
\begin{equation}
\frac{\alpha(t)}{t}\ \mbox{is a non-decreasing function}
\end{equation}
and let $\lambda:\mathbb{N}\to\mathbb{N}$ be a function. Let $\mathcal{Q}$ be a collection of pairwise disjoint cubes $Q_i\subset \mathbb{R}^n\setminus E$. We define for each such collection $\mathcal{Q}$ and for every $k\in \mathbb{N}$ a function
$$\chi_k^{\mathcal{Q}}(x)= \left\{ \begin{array}{l}
    1, \ \mbox{if one can find cubes}\ Q_i^k(x)\in \mathcal{Q},\  i=1,...,\lambda(k),\  \\
\mbox{such that } Q_i^k(x)\subset A_k(x)\ \mbox{and}\  \diam(Q_i^k(x))\ge \alpha(2^{-k})\  \mbox{for all } i;\\
    0, \ \mbox{otherwise}.  \end{array} \right.$$
Here $A_k(x)=\{y\in \mathbb{R}^n:2^{-k}<|x-y|<2^{-k+1}\}$. Let
$$S_j^{\mathcal{Q}}(x) = \sum_{k=1}^j\chi_k^{\mathcal{Q}}(x).$$
We say that a set $E$ is weakly mean porous with parameters $\alpha$ and $\lambda$,
if there exists a collection $\mathcal{Q}$ as above and an integer $j_0\in \mathbb{N}$
such that
\begin{equation}\label{porositycondition}
\frac{S_j^{\mathcal{Q}}(x)}{j}\ >\ \frac{1}{2}
\end{equation}
for all $x\in E$ and for all $j\ge j_0$.
\end{definition}

Let us remark that the concepts of porosity and mean porosity are well-known tools in geometric analysis, see e.g. \cite{KR}. In the proof of our main result, we will apply the following sharp estimate on the Hausdorff (gauge) dimension of weakly mean porous sets established in \cite[Corollary 3.5]{N}.
\begin{lemma}[\cite{N}]\label{poroustheorem}
Let $E\subset \mathbb{R}^n$ be a weakly mean porous set with parameters $\alpha$ and $\lambda$
such that
\begin{equation}\label{porousmonotonity}
\frac{\lambda(k)\alpha(2^{-k})^n}{(2^{-k})^n}\ \mbox{is a non-increasing function of}\  k
\end{equation}
and
\begin{equation}\label{porousdivergence}
\sum_{k=j_0}^{\infty}\frac{\lambda(k)\alpha(2^{-k})^n}{(2^{-k})^n}=\infty.
\end{equation}
Then $m_n(E)=0$ and, moreover, there is a positive constant $C(n)$ such that $\Lambda_{h}(E)=0$ for each premeasure $h$, which satisfies
$$h(2^{-j})\le M2^{-jn}\exp\Big(C(n)\sum_{k=j_0}^j\frac{\lambda(k)\alpha(2^{-k})^n}{(2^{-k})^n}\Big)$$
for all $j>j_0$ with some positive constant $M$.
\end{lemma}

Throughout the proofs we denote by $C$ positive constants depending only on the given data $n,K$ and $\beta$. These constants may vary from expression to expression as usual.

The \emph{average derivative} $a_f$ of a quasiconformal mapping $f$ is defined by
$$a_f(x)=\Big(\frac{1}{m_n(B_x)}\int_{B_x} J_f\ dm_n\Big)^{1/n},$$
where $x\in \B^n$ and $B_x=B(x,\frac{1}{2}(1-|x|))$. Recall that $a_f$ satisfies a Harnack inequality
\begin{equation}\label{harnack}
1/C\le \frac{a_f(z)}{a_f(y)}\le C
\end{equation}
for any points $z,y$ belonging to some Whitney-ball $B_x\subset \B^n$, see e.g. \cite{BKR}. Also note that
\begin{equation}\label{afdf}
1/C\int_Q a_f^n\ dm_n\le \int_Q|Df|^n dm_n\le C\int_Q a_f^n\ dm_n
\end{equation}
for all cubes $Q$ in a Whitney decomposition $\mathcal{W}$ of $\B^n$, see e.g. \cite[Theorem 3.4]{AK}. A Whitney decomposition of $\B^n$ refers to a collection of closed dyadic cubes $Q\subset \B^n$ with pairwise disjoint interiors such that
$$\bigcup_{Q\in \mathcal{W}}Q=\B^n$$
and that $\diam(Q)\le \dist(Q,\partial \B^n)\le 4\diam(Q)$. See \cite{S} for the existence of such a decomposition.

Let us also recall that the \emph{quasihyperbolic distance} $k_{\Omega}(x_1,x_2)$ between two points $x_1,x_2$ in a domain $\Omega\subsetneq \mathbb{R}^n$ is defined as the infimum of
$$\int_{\gamma}\frac{ds}{\dist(x,\partial \Omega)}$$
over all rectifiable curves joining $x_1$ to $x_2$ in $\Omega$.

Before the proof of our main result, we prove an additional technical property for the weight function $\varphi$.
\begin{lemma}\label{tech}
Let $\varphi(t)$ be a weight function satisfying the conditions \eqref{prop1}, \eqref{prop2}, \eqref{prop3}. Then there exists a constant $C=C(n,\beta)>0$ such that
$$\int_{[0,r]}\frac{\varphi(t)^{1/n}}{t}dt \le \varphi(r)^C$$
for all sufficiently small $r>0$.
\end{lemma}

{\it Proof.} Let us write $\varphi(t)=\exp(-\alpha(t))$. We show first that
\begin{equation}\label{first}
\log\log\frac{1}{r} = o(\alpha(r))\ \ \ \ \mbox{as}\ r\to 0.
\end{equation}
Suppose that this assertion fails. Then there is a constant $c>0$ such that for an arbitrarily small $r>0$ we have that
$$\varphi(r)\ge \exp(-c\log\log\frac{1}{r})=(\log\frac{1}{r})^{-c}.$$
This implies that
$$\varphi^{-1}(s)\le \exp(-s^{-1/c})$$
or equivalently
$$\log\frac{1}{\varphi^{-1}(s)}\ge s^{-1/c}$$
for some arbitrarily small $s>0$.
This, however, contradicts the condition \eqref{prop3}. Thus we have proven \eqref{first}. 

On the other hand, it follows from the assumptions with the help of Gronwall's lemma \cite[p. 436]{W} that, for all sufficiently small $t>0$,
$$\varepsilon \le \Big|\frac{\alpha'(t)t\log t}{\alpha(t)}\Big|\le \frac{1}{\varepsilon}$$
with some $\varepsilon>0$ depending only on $\beta$ (cf. the proof of \cite[Remark 5.3]{N}). By combining this with \eqref{first} we obtain
\begin{align}
\frac{\varphi(t)^{1/n}}{t}=\frac{\exp(-\frac{1}{n}\alpha(t))}{t} &\le \frac{\exp(-\frac{1}{2n}\alpha(t)+\log(\alpha(t))-\log\log\frac{1}{t})}{t} \notag \\
&=\frac{\alpha(t)}{t\log\frac{1}{t}}\exp(-\frac{1}{2n}\alpha(t)) \notag\\
&\le C|\alpha'(t)|\exp(-\frac{1}{2n}\alpha(t)) \notag
\end{align}
for all sufficiently small $t>0$. The claim follows.

{\it Proof of Theorem \ref{main}}. Let $E_{\infty}$ consist of those points $\xi\in \partial \B^n$ for which $\int_{[0,\xi]}a_f(x)\ |dx|=\infty$. Then $\capacity_n(E_{\infty})=0$ by \cite[Theorem 4.4]{BKR} and thus, by Lemma \ref{tech}, $\Lambda_{\varphi}(E_{\infty})=0$, see e.g. \cite[p. 120]{Re} or \cite[Remark 1.3]{NT}. Let us then write
$$G_j=\{\xi\in \partial \B^n\setminus E_{\infty}:\int_{[0,\xi]}a_f(x)\ |dx|\le j\}$$
for $j\in \mathbb{N}$, and let $S_j$ consist of the union of Stolz cones at $G_j$, i.e.,
$$S_j=\bigcup_{\xi\in G_j}\bigcup_{0\le t<1}B(t\xi,\frac{1}{2}(1-t)).$$
Then $S_j$ is open and $\diam f(S_j)\le C\sup_{\xi\in G_j}\{\int_{[0,\xi]}a_f(x)\ |dx|\}\le Cj$, where the first inequality is implied by the Harnack inequality \eqref{harnack}, \cite[Remark 3.3]{BKR}, and the fact that the discrete length of a curve (as defined in \cite{HR}) is always at least the Euclidean distance between the end points. Thus we have that $m_n(f(S_j))<\infty$. It follows by the inequalities \eqref{distortion} and \eqref{afdf} that the function $u_j:\B^n\to (0,\infty)$,
$$u_j(x)= \left\{ \begin{array}{cl}
    a_f(x)^n & \mbox{for}\ x\in S_j;\\
    0 &  \mbox{otherwise}  \end{array} \right.$$
belongs to $L^1(\B^n)$. This implies with the help of Besicovitch's covering theorem that there is a set $E_j\subset G_j$ with $\Lambda_{\varphi}(E_j)=0$ so that
$$\int_{B(\xi,r)\cap \B^n} u_j\  dm_n = o(\varphi(r))\ \mbox{as}\ r\to 0$$
for all $\xi\in G_j\setminus E_j$ (cf. \cite[p. 118]{Z}). In particular, for all $\xi \in G_j\setminus E_j$,
$$\int_{B_{t\xi}}a_f^n\ dm_n=\int_{B_{t\xi}}u_j\ dm_n\le \int_{B(\xi,\frac{3}{2}(1-t))}u_j\ dm_n=o(\varphi(1-t))\ \mbox{as}\ t\to 1$$
and hence, by the Harnack inequality \eqref{harnack},
\begin{equation}\label{linfty}
a_f(t\xi)^n\le C(1-t)^{-n}\int_{B_{t\xi}}a_f^n\ dm_n=o\Big(\frac{\varphi(1-t)}{(1-t)^n}\Big)\ \mbox{as}\ t\to 1.
\end{equation}

Let us define $E=E_{\infty}\cup \bigcup_j E_j$. Then $\Lambda_{\varphi}(E)=0$ by the subadditivity of $\varphi$-measure. We claim that $f(\partial \B^n\setminus E)$ is of zero $\psi$-measure. Since $\partial \B^n\setminus E\subset \bigcup_j (G_j\setminus E_j)$, it suffices to show that, for all $j\in \mathbb{N}$, $f(G_j\setminus E_j)$ is of zero $\psi$-measure.

For each $m\in \mathbb{N}$ we define a set $F_j^m\subset G_j\setminus E_j$  by
$$F_j^m=\{\xi\in G_j\setminus E_j: a_f(t\xi)\le \frac{\varphi(1-t)^{1/n}}{1-t}\ \mbox{for all}\ t\ge 1-2^{-m}\},$$
whence $\bigcup_m F_j^m=G_j\setminus E_j$ by the inequality \eqref{linfty}. Moreover, $F_j^1\subset F_j^2\subset ...$, and thus, by the subadditivity of the Hausdorff measure, it suffices to show that $f(F_j^m)$ is of zero $\psi$-measure for an arbitrarily large integer $m$ in order to prove the theorem.

Because of Lemma \ref{poroustheorem} it now suffices to prove that the set $f(F_j^m)\subset \partial \Omega$ is weakly mean porous with parameters $C\alpha$ and $C\lambda$, where (for small $t$)
\begin{equation}\label{parameters}
\alpha(t)=cu(t)/u'(t)\ \ \mbox{and}\ \ \lambda(k)\ge 2^{-k}/\alpha(2^{-k})
\end{equation}
and $c>0$ depends only on $n,K$ and $\beta$.

Let $j,m\in \mathbb{N}$ and let $\xi\in F_j^m$. Then the radial limit $f(\xi)$ exists by \cite[Remark 4.5]{BKR} and, by the Harnack inequality \eqref{harnack} and Lemma \ref{tech}, we obtain
\begin{align}\label{modcont}
|f(t\xi)-f(\xi)| &\le \sum_{Q\in \mathcal{W}: Q\cap [t\xi,\xi]\neq \emptyset} \diam f(Q) \le C\int_{[t\xi,\xi]}a_f(x)\ |dx| \notag\\
&\le C\int_{[t,1]}\frac{\varphi(1-s)^{1/n}}{1-s}ds \le \varphi(1-t)^C
\end{align}
for all $t\ge 1-2^{-m}$, provided that $m$ was chosen large enough above. Here we denoted by $\mathcal{W}$ a Whitney decomposition of $\mathbb{B}^n$.

Let us write $\gamma=f([0,\xi))$. We choose an integer $i_0\ge m$ so large that $2^{-i_0+1}\le \dist(f(0),\partial f(\B^n))$, and we define for all integers $k\ge i_0$ a function
$$\chi_k(f(\xi))= \left\{ \begin{array}{l}
    1, \ \mbox{if } k_{\Omega}(f(t_a\xi),f(t_b\xi))\le
C\frac{2^{-k}}{\alpha(2^{-k})}\\
    0, \ \mbox{otherwise},  \end{array} \right.$$
where $f(t_a\xi)$ and $f(t_b\xi)$ are the last entry point along $f([0,\xi])$ into $A_k(f(\xi))$ and the first exit point after $f(t_a\xi)$ along $f([0,\xi])$ from $A_k(f(\xi))$, respectively, and $k_{\Omega}$ is the quasihyperbolic metric in $\Omega$.
We also define for all integers $i\ge i_0$ a function
$$S_i(f(\xi))=\sum_{k=i_0}^i\chi_k(f(\xi)).$$
We then claim that
\begin{equation}\label{sjclaim}
\frac{S_i(f(\xi))}{i}>\frac{1}{2}
\end{equation}
for all sufficiently large $i\in \mathbb{N}$ provided that $c$ in \eqref{parameters} is chosen small enough.

Let us consider an annulus $A_k(f(\xi))$ such that $\chi_k(f(\xi))=0$. The quasihyperbolic distance $k_{\Omega}$ between the points $f(t_a\xi)$ and $f(t_b\xi)$ is at least $C2^{-k}/\alpha(2^{-k})$. Due to the quasi-invariance of the quasihyperbolic metric under quasiconformal mappings \cite[p. 62]{GO}, we then have that the quasihyperbolic distance $k_{\B^n}(t_a\xi,t_b\xi)$ is at least $C2^{-k}/\alpha(2^{-k})$ (with $C$ still depending only on $n,K$ and $\beta$).

Consider the largest $t<1$ with
$$|f(t\xi)-f(\xi)|=2^{-i}.$$
It follows from \eqref{modcont} that $2^{-i}\le \varphi(1-t)^C$ and hence, by \eqref{prop3},
\begin{equation}\label{yla}
\log\frac{1}{1-t}\le\log\Big(\frac{1}{u(2^{-i/C})}\Big) \le C\log\Big(\frac{1}{u(2^{-i})}\Big).
\end{equation}
On the other hand, we have by the observations above that
\begin{equation}\label{ala}
\log\frac{1}{1-t}=k_{\B^n}(0,t\xi)\ge \sum k_{\B^n}(t_a\xi,t_b\xi)\ge \sum C\frac{2^{-k}}{\alpha(2^{-k})},
\end{equation}
where the summation is over all $i_0\le k\le i$ with $\chi_k(f(\xi))=0$.

Suppose that the assertion \eqref{sjclaim} fails for some large integer $i$. Then, by combining \eqref{yla} and \eqref{ala}, we arrive at (the summation indices follow from the assumption \eqref{prop2})
\begin{align}
\log\Big(\frac{1}{u(2^{-i})}\Big) &\ge C\sum_{k=i_0}^{i/2}\frac{2^{-k}u'(2^{-k})}{cu(2^{-k})} \notag\\
&\ge C\frac{1}{c}\Big(\log\Big(\frac{1}{u(2^{-i/2})}\Big)-\log\Big(\frac{1}{u(2^{-i_0})}\Big)\Big). \notag
\end{align}
But this inequality is a contradiction with property \eqref{prop3} if we choose $i$ large enough and $c$ small enough (depending on $\beta$). This proves \eqref{sjclaim}.

For the final arguments of the proof we define a collection $\mathcal{Q}$ of disjoint cubes in the domain $\Omega$ in the following way. Let $\mathcal{W}$ be a Whitney decomposition of $\Omega$. Then let $\mathcal{Q}$ consist of the interiors of all the cubes in the Whitney decompositions of the cubes $Q\in \mathcal{W}$. We claim that
\begin{equation}\label{chiclaim}
\chi_k^{\mathcal{Q}}(f(\xi)) \ge \chi_k(f(\xi)).
\end{equation}
for all $k\ge i_0$ with parameters $C\alpha$ and $C\lambda$ (defined above).

Let us consider a ``good'' annulus, i.e., $A_k(f(\xi))$ with $\chi_k(f(\xi))=1$. Then
$$k_{\Omega}(f(t_a\xi),f(t_b\xi))\le C\frac{2^{-k}}{\alpha(2^{-k})},$$
which means geometrically that there are at most $C2^{-k}/\alpha(2^{-k})$ Whitney cubes $Q\in \mathcal{W}$ intersecting the quasihyperbolic geodesic joining $f(t_a\xi)$ and $f(t_b\xi)$. But since the length of this curve is at least the width of the annulus $A_k(f(\xi))$ or $2^{-k}$, it follows that some of these Whitney cubes must have a large diameter. Indeed, one finds at least $C2^{-k}/\alpha(2^{-k})$ cubes $Q\in \mathcal{Q}$ so that $Q\subset A_k(f(\xi))$ and the diameter of each cube is at least $C\alpha(2^{-k})$. This observation follows by easy geometric arguments and elementary calculations, which we leave to the reader (cf. \cite[Lemma 4.6]{N}). Let us point out that the second Whitney decomposition $\mathcal{Q}$ is needed here to ensure that there are enough cubes in $A_k(f(\xi))$ also in the case that the geodesic intersects only, say, one large Whitney cube from the first decomposition $\mathcal{W}$.

In conclusion, we have shown \eqref{chiclaim} and \eqref{sjclaim}, which together imply that the set $f(F_j^m)$ is weakly mean porous. The desired estimate on the $\psi$-measure of $f(F_j^m)$ (and thus for the $\psi$-measure of $f(\partial \B^n\setminus E)$) then follows from Lemma \ref{poroustheorem}, and thus the proof is complete.

\section{Sharpness of the results}\label{sharpnesssection}

In this section we show the sharpness of Theorem \ref{main} in the interesting case that $\Lambda_{\varphi}$ is absolutely continuous wit respect to $H^1$. At the same time we show the ``only if''-part of Theorem \ref{maincorollary}. The example domain that we consider here was constructed also in \cite{N} for a different purpose. This example turns out to be the critical one also for the results of this paper. We give an outline for the construction in the following. See \cite[Section 7.2]{N} for the full details.

Given a gauge function $\varphi$ satisfying \eqref{prop1}, \eqref{prop2} and \eqref{prop3} we choose an increasing function $\alpha$ so that
$$\alpha(2^{-k})=cu(2^{-k})/u'(2^{-k})$$
for all $k\in \mathbb{N}$, where the constant $c>0$ is to be determined below. We assume (by \eqref{prop2}) that $\alpha(t)\le t/16$ for all $t>0$. Moreover, we assume that $\alpha(2^{-k})$ is dyadic for all $k$.

Starting with the square $Q_1=\{(x,y)\in \mathbb{R}^2:|x|< 2^{-1}\ \mbox{and}\ |y|<2^{-1}\}$ define $\Omega_1$ as the intersection of $Q_1$ and the open $\alpha(2^{-1})$-neighborhood of the coordinate axes in $Q_1$. Then subdivide $Q_1$ into $4$ dyadic squares $Q_2^i$, $i\in \{1,2,3,4\}$, and define $\Omega_2$ as the union of the intersections of $Q_1\setminus \overline{\Omega}_1$ and the open $\alpha(2^{-2})$-neighborhoods of the centered coordinate axes of each square $Q_2^i$. Then attach each component of $\Omega_2$ to $\Omega_1$ in the way shown in the picture below.

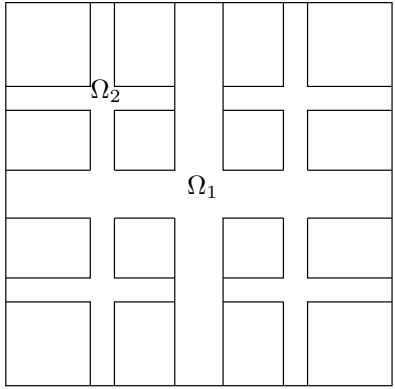
\begin{center}
\begin{picture}(0,0)%
\includegraphics{cantor.pstex}%
\end{picture}%
\setlength{\unitlength}{3947sp}%
\begingroup\makeatletter\ifx\SetFigFont\undefined%
\gdef\SetFigFont#1#2#3#4#5{%
  \reset@font\fontsize{#1}{#2pt}%
  \fontfamily{#3}\fontseries{#4}\fontshape{#5}%
  \selectfont}%
\fi\endgroup%
\begin{picture}(2424,2424)(589,-2173)
\put(1726,-961){\makebox(0,0)[lb]{\smash{{\SetFigFont{10}{12.0}{\rmdefault}{\mddefault}{\updefault}{\color[rgb]{0,0,0}$\Omega_1$}%
}}}}
\put(1126,-361){\makebox(0,0)[lb]{\smash{{\SetFigFont{10}{12.0}{\rmdefault}{\mddefault}{\updefault}{\color[rgb]{0,0,0}$\Omega_2$}%
}}}}
\end{picture}%

\end{center}

Continue the construction by subdividing each square $Q_2^i$ into $4$ dyadic squares and defining $\Omega_3$ accordingly. By iterating this process one obtains a simply connected domain $\Omega$ (take $\bigcup_j \Omega_j$ and open certain gates as in the picture above to make it simply connected), which satisfies the growth condition
$$k_{\Omega}(0,x)\le \phi\Big(\frac{\dist(x,\partial \Omega)}{\dist(0,\partial \Omega)}\Big) + C$$
on the quasihyperbolic metric $k_{\Omega}$ with the function $\phi(t)=\frac{C}{c}\log\frac{1}{u(t)}$. Then by Lemma \ref{tech}, we can use known results from the literature (\cite[Theorem 1.2]{HK}) to conclude that there exists a (quasi)conformal mapping $f:\B^2\to \Omega$, which is uniformly continuous with a modulus of continuity $C\varphi(t)$, i.e.,
\begin{equation}\label{examplemodulus}
|f(x)-f(y)|\le C\varphi(|x-y|)
\end{equation}
for all $x,y\in \B^2$ provided that the constant $c$ is chosen large enough in the construction above.

Moreover, by standard arguments (cf. \cite[Section 7.2]{N} or \cite{KR}) involving a construction of a ``Frostman measure'' on the boundary of $\Omega$ and an employment of the Frostman's lemma one observes that $\Lambda_{\psi}(\partial \Omega)>0$ with a dimension gauge $\psi(r)$ satisfying
$$\psi(r)\le r^2\exp\Big(\tilde{C}\int_{[r,r_0]}\frac{u(t)}{u'(t)}\frac{dt}{t^2}\Big)$$
for all $r>0$ with a sufficiently large constant $\tilde{C}>0$. In particular, this means that $m_2(\partial \Omega)>0$ if the integral condition \eqref{divergence} fails.

It only remains to show that the $\psi$-measure of the set $f(\partial \B^2\setminus E)$ is also positive for any set $E\subset \partial \B^2$ of zero $\varphi$-measure. To that end, let $E\subset \partial \B^2$ be an arbitrary set with $\Lambda_{\varphi}(E)=0$. This means, by definition, that for any $\varepsilon>0$ there is a countable collection of balls $B_1,B_2,...$ with diameters $r_1,r_2,...$ so that $E\subset \bigcup_i B_i$ and 
\begin{equation}\label{measureE}
\sum_{i}\varphi(r_i)< \varepsilon.
\end{equation}

On the other hand, note that the internal diameter of $\Omega$ is finite (recall that the internal distance between two points $z,w \in \Omega$ is the infimum of lengths of curves $\gamma :[0,1] \to \Omega$ such that $\gamma (0)=z$ and $\gamma (1)=w$), and hence the radial limit $f(\xi)$ exists for all points $\xi\in \partial \B^2$. This follows from the Gehring--Hayman theorem, see e.g. \cite[Remark 4.5]{BKR}. Thus $f(E)\subset \partial \Omega$ is well defined, and we observe that $f(E)\subset \bigcup_i f(B_i)$. Furthermore, $\diam(f(B_i))\le C\varphi(r_i)$ for each ball $B_i$ by \eqref{examplemodulus}. By combining this with \eqref{measureE} we conclude that $H^1(f(E))=0$ which implies $\Lambda_{\psi}(f(E))=0$, since we are assuming that $\Lambda_{\psi}$ is absolutely continuous with respect to $H^1$, which as we said above, is the interesting case. It follows that $\Lambda_{\psi}(f(\partial \B^2\setminus E))\ge \Lambda_{\psi}(\partial \Omega\setminus f(E))>0$.

Let us finally point out that the construction of $\Omega$ above can be extended to $\mathbb{R}^n$, $n\ge 3 $ so that the resulting domain is quasiconformally equivalent to the ball $\B^n$, see \cite{V}. For example, in the case $n=3$ define $\Omega_1$ in the unit cube $Q_1$ (of side length $1$) by taking the $\alpha(2^{-1})$-neighborhood of the coordinate axes and of the lines $(t,\pm 2^{-2},0)$ (we have to include the neighborhoods of these additional lines to make the final set $\Omega$ connected). Then subdivide $Q_1$ into $8$ dyadic (open) subcubes $Q_2^i$ and define $\Omega_2\subset \bigcup Q_2^i\setminus \overline{\Omega}_1$ accordingly. Iterate this process and, in the end, attach each component of $\Omega_{j+1}$ to $\Omega_{j}$ to obtain a connected set $\Omega$. One can then show the sharpness of Theorem \ref{main} with similar calculations and conclusions as in the planar case discussed above.

{\it Acknowledgement.} The authors wish to thank Professor Pekka Koskela for pointing out the problem studied in this paper.

\smallbreak

Tomi Nieminen, Department of Mathematics and Statistics, P.O. Box 35, FIN-40014 University of Jyv\"askyl\"a, Finland.\\ \textit{E-mail address:} \ \texttt{tominiem@maths.jyu.fi}

Ignacio Uriarte-Tuero, Mathematics Department, 202 Mathematical Sciences Bldg., University of Missouri, Columbia, MO 65211-4100, USA; and\\ Fields Institute, 222 College Street, Toronto, Ontario, M5T 3J1, Canada.\\ \textit{E-mail address:} \ \texttt{ignacio@math.missouri.edu}

\end{document}